\documentclass[a4paper]{article}
\long\def\kmcomment#1{}
\usepackage{graphicx,color,verbatim,setspace,ntheorem} 
\usepackage{latexsym,amssymb,amsmath}
\usepackage[paperwidth=9.0in, paperheight=10.5in]{geometry}
\newcommand{\myd}{d\, }

\newcommand{\kmqed}{\hfill \rule{1ex}{1.5ex}\par}

\newcommand{\Shat}[1]{ \text{S}_{#1}}
\newcommand{\cb}{\text{cbase}}
\newcommand{\ol}[1]{$\ensuremath{\overline{#1}}$}
\newcommand{\ool}[1]{\overline{#1}}
 
\newcommand{\pos}[2]{\text{pos}_{#1}(#2)} 
\newcommand{\CCx}[2]{\ensuremath{\mathfrak{C}^{#1}_{#2}}} 
\newcommand{\ds}{\displaystyle }
\newcommand{\frakS}[1]{\mathfrak{S}_{#1}}
\newcommand{\fraksp}{\mathfrak{s}\mathfrak{p}}
\newcommand{\frakham}{\mathfrak{h}\mathfrak{a}\mathfrak{m}}
\newcommand{\Pkt}[2]{\{#1,#2\}}%
\newcommand{\mR}{\ensuremath{\mathbb{R}}} 
\newcommand{\mN}{\ensuremath{\mathbb{N}}} 
\newcommand{\rank}{\text{rank}} 
\newcommand{\red}{black} 

\newcommand{\HGFs}[2]{\text{H}^{#1}_{\rm GF}({\frakham}_6^0,{Sp}(6,\mR))_{#2}}

\newcommand{\CGFF}[4]{\text{C}^{#1}_{GF}({\frakham}_{#2}^{#3})_{#4}} 
\newcommand{\HGFF}[4]{\text{H}^{#1}_{GF}({\frakham}_{#2}^{#3})_{#4}}

{
\theorembodyfont{\rmfamily}
\theoremstyle{plain}
\newtheorem{thm}{Theorem}[section] 

\newtheorem{exam}{Example}[section]
\newtheorem{defn}{Definition}[section]
\newtheorem{kmProp}{Proposition}[section]
\newtheorem{Lemma}{Lemma}[section]
 
\newtheorem{kmCor}[kmProp]{Corollary}
\newtheorem{kmRemark}{\textbf{Remark}}[section]
}

\newcommand{\CGF}[4]{\text{C}^{#1}_{\rm GF} ( {\frakham}_{#2}^{#3}, {\fraksp}(#2,\mR))_{#4}}
\newcommand{\HGF}[4]{\text{H}^{#1}_{\rm GF} ( {\frakham}_{#2}^{#3}, {\fraksp}(#2,\mR))_{#4}}
\newcommand{\cgf}[4]{\text{C}^{#1}_{\rm GF} ( {\frakham}_{#2}^{#3})_{#4}}

\title{The relative 
Gel'fand-Kalinin-Fuks cohomology groups of the formal Hamiltonian
vector fields on 6-dimensional plane}

\author{Kentaro Mikami\thanks{%
  Akita University, supported by Grant-in-Aid for 
  Scientific Research (C) of JSPS, No.23540067 and No.20540059}
} 
\date{\today} 

\begin{document}
\onehalfspacing
{\allowdisplaybreaks 
\maketitle

\section{Introduction}%
At the beginning of research of 
the (relative) Gel'fand-Fuks cohomology group,  
Gel'fand-Kalinin-Fuks (\cite{MR0312531}) got that 
$\displaystyle \HGF{\bullet}{2}{ }{w} = 0$ for $w=2,4,6$ and the 
$\displaystyle \HGF{7}{2}{ }{8} \cong \mR$, whose generator is called    
the Gel'fand-Kalinin-Fuks class.   
The next non-trivial result in this context is 
$\displaystyle \HGF{9}{2}{ }{14} \cong \mR$, which is discovered by 
S.~Metoki (\cite{metoki:shinya}) in 1999.      
In short, $\displaystyle \frakham_{2}^{ }$ is the Lie algebra of the
Hamiltonian vector fields of the formal polynomials on 
$\displaystyle \mR^2$.  

D.~Kotschick and S.~Morita (\cite{KOT:MORITA}) 
research 
$\displaystyle \HGF{\bullet}{2}{0}{w} $ and determined the whole spaces
while $w\le 10$, where   
$\displaystyle \frakham_{2}^{0}$ is the Lie subalgebra of the
Hamiltonian vector fields of the formal polynomials which vanish at the origin
of $\displaystyle \mR^2$.  

Inspired by  \cite{KOT:MORITA}, we are interested in studying higher weight or
higher dimensions.  When $n=1$, we have got some results for higher weight
cases in \cite{M:N:K}, and when $n=2$ for lower weight cases in
\cite{Mik:Nak}.  
A big difference of methodology in the two cases happens when getting the
irreducible decomposition of the tensor product of irreducible representations
of $Sp(2 n,\mR)$.  When $n=1$, it is lucky we can use the Clebsch-Gordan
rule.  When $n=2$, in \cite{Mik:Nak} we have used the Littlewood-Richardson rule. 

In this paper, we deal with the case of $n=3$.  
The main target is again the space of homogeneous polynomials.  
If we put 
$\displaystyle \Shat{k}(\mR^{2n})$ be the vector space of order $k$ homogeneous polynomials
of $\displaystyle x_1,x_2, \ldots, x_{2n}$ of $\displaystyle \mR^{2n}$,    
then it is known that 
$\displaystyle \dim \Shat{k}(\mR^{2n}) = \frac{ (2n-1+k)!} { (2n-1)! k! }$.  
The table below shows how fast the dimension of $\displaystyle
\Shat{k}(\mR^{2n})$ increases in accordance with $k$  or $n$:  
\begin{center}
\begin{tabular}{c| *{7}{r}}
        k & 1 & 2 & 3 & 4 & 5 & 6 & $\cdots$ \\\hline
        $n=1$ & 2 & 3 & 4 & 5 & 6 & 7 &$\cdots$ \\\hline
        $n=2$ & 4 & 10 & 20 & 35 & 56 & 84& $\cdots$ \\\hline
        $n=3$ & 6 & 21 & 56 & 126 & 252 & 462& $\cdots$ 
\end{tabular}
\end{center}

By the methodology we have used so far, we encountered some difficulty in
computations when $n=3$.  So,  this time, instead of the Littlewood-Richardson
rule, we would like to use of the crystal base theory by M.~Kashiwara
(\cite{kashi:french}, \cite{kashi:nakashima}, \cite{nakashima}).  

\bigskip

Abbreviating  
 the relative cochain space   
 $\ds \CGF{j}{6}{0}{w}$ by $\ds \CCx{j}{w}$, we introduce our first result and 
the second result about Betti numbers for weight 2, 4 and 6.  
 
\textbf{Theorem 1:} 
Dimensions of each relative cochain complexes of weight 2 or 4 are as below:
\begin{center}
        \begin{tabular}{c|*{2}c}
                & $\ds\CCx{1}{2}$ & $\CCx{2}{2}$ \\\hline
                $\dim$ & 0 & 1
        \end{tabular}
\qquad and 
\qquad 
        \begin{tabular}{c|*{4}c}
                & $\ds\CCx{1}{4}$ & $\CCx{2}{4}$ & $\ds\CCx{3}{4}$ & $\CCx{4}{4}$
                \\\hline
                $\dim$ & 0 & 0& 1 & 3
        \end{tabular} 
\end{center}
The Euler characteristic number of $\ds \HGFs{\bullet}{2} = 1$ and  
the Euler characteristic number of $\ds \HGFs{\bullet}{4} = 2$.   

\medskip

The table of dimensions of each relative cochain complexes with weight 6 is
\begin{center}
        \begin{tabular}{c|*{6}c}
                & $\ds\CCx{1}{6}$ & $\CCx{2}{6}$ 
                & $\ds\CCx{3}{6}$ & $\CCx{4}{6}$
                & $\ds\CCx{5}{6}$ & $\CCx{6}{6}$
                \\\hline
                $\dim$ & 0 & 1& 1 & 0& 4 & \textbf{7}  
        \end{tabular}
\end{center}
and 
the Euler characteristic number of $\ds \HGFs{\bullet}{6} = 3$.   

\medskip

The Euler characteristic number means in this paper  the alternating sum of
Betti numbers except 0-dimensional. 

\bigskip


\textbf{Theorem 2:} The Betti numbers of weight 2 relative
Gel'fand-Kalinin-Fuks cohomology groups are $\ds b^{0}_{2}=1,\ b^{1}_{2}=0,\
b^{2}_{2}=1$. 

About the Betti numbers of  weight 4 relative Gel'fand-Kalinin-Fuks cohomology
groups, $\ds b^{4}_{4} = 2,\ b^{0}_{4}=1$ and the others are 0. 

In the case of weight 6, we have
$\ds b^{0}_{1}=1,\ b^{1}_{6} = b^{2}_{6} = b^{3}_{6}=b^{4}_{6}=0,\  
b^{5}_{6}=2,\ \text{and} \ b^{6}_{6} = 5$.

\section{Preliminary}
We are
interested in the relative Gel'fand-Kalinin-Fuks cohomology groups of 
$\displaystyle \frakham_{2n}^{0}$ when $n=3$.  
We review the notions we deal with in this paper, but we refer the
precise definitions to \cite{M:N:K}.   

We can split the polynomial functions by
their homogeneity, and we regard the Lie algebras as follows.  
\begin{align*}
        \ds  \frakham_{2n}^{ } =& \left(\mathop{\oplus}_{p=1}^{\infty}
\Shat{p}(\mR^{2n}) \right)^{\wedge} \quad \text{is a Lie algebra }\\
\frakham_{2n}^{0} =& \left(\mathop{\oplus}_{p=2}^{\infty}
\Shat{p}(\mR^{2n})\right)^{\wedge} \quad \text{is a subalgebra of } 
\frakham_{2n}^{ } \\ 
\noalign{and}
\frakham_{2n}^{1} =& \left(\mathop{\oplus}_{p=3}^{\infty}
\Shat{p}(\mR^{n2})\right)^{\wedge} \quad \text{is a subalgebra of } 
\frakham_{2n}^{ }  
\end{align*} 
where $\ds \left(\phantom{MM}\right)^{\wedge}$ means the completion with the
Krull topology.

The cochain complex is the exterior algebra of
the dual spaces of $\ds\Shat{p}(\mR^{2n})$'s, and we have the notion of ``weight'' on the cochain
complex. 

\begin{defn} 
{\color{black}Let} $\displaystyle \frakS{\ell}$ be the dual space of
$\ell$-homogeneous polynomial functions $\displaystyle \Shat{\ell}(\mR^{2n})$,
and define the weight of each non zero element of
$\displaystyle \frakS{\ell}$ to be $\ell-2$.  
For each non-zero element of $\displaystyle 
\frakS{\ell_1} \wedge 
\frakS{\ell_2} \wedge \cdots \wedge 
\frakS{\ell_s}$ 
($ \ell_1 \le \ell_2 \le \cdots \le \ell_s$), 
we 
define  its weight to be $\displaystyle \sum_{i=1}^{s} (\ell_i -2)$. 
\end{defn}

\begin{kmProp}[cf.\cite{KOT:MORITA},\cite{M:N:K}]
The coboundary operator $\myd$ of the Gel'fand-Kalinin-Fuks cochain
complex preserves the weight,
namely, 
if a cochain $\sigma$ is of weight $w$, then 
$\displaystyle \myd \sigma$ is also of weight $w$.

Hence we can decompose the total space of cochain complex by degree and
weight: namely, 
\begin{align*}
\CGFF{m}{2n}{ }{w}
= & \text{Linear Span of }\{  
\sigma \in 
\Lambda^{k_1} 
\frakS{1} \wedge 
\Lambda^{k_2} 
\frakS{2} \wedge \cdots
\wedge 
\Lambda^{k_s} 
\frakS{s} \\ & \hspace{35mm}
\mid 
\sum_{j=1}^{s} k_j = m \ , \  
\sum_{j=1}^{s} k_j (j-2) = w{\color{black}\ , \ s=1,2,\ldots} 
\}\end{align*} 
and we can define the cohomology group 
$ \displaystyle 
\HGFF{m}{2n}{ }{w}
$.  
\end{kmProp}

$\displaystyle 
\CGFF{\bullet}{2n}{0}{w}
$ is the subspace of    
$\displaystyle 
\CGFF{\bullet}{2n}{ }{w}
$ characterized by 
$k_1=0$.
If we restrict our attention to the cochain complex relative to
$Sp(2n,\mR)$,
then it turns out $k_2=0$ 
and 
\begin{equation}
        \label{eqn:Kont:} 
        \ds\CGF{\bullet}{2n}{0}{w} =\left( \cgf{\bullet}{2n}{1}{w}
\right)^{triv}
\end{equation} 
(cf.\cite{M:N:K}){\color{black}.}

\begin{kmRemark}[cf.\cite{KOT:MORITA}]
Since we have the negative Identity matrix in $\ds Sp(2n,{\color{black}\mR})$, we see that the relative cochain
complex of odd weight must be the zero space, and hence we only deal
with the complexes of even weights.  
\end{kmRemark}

From (\ref{eqn:Kont:}), when we study $\ds \CGF{\bullet}{2n}{0}{w}$,   
we first look at the subcomplex
$\ds\cgf{\bullet}{2n}{1}{w}$, namely, $\ds k_1=k_2=0$.  So,  
for $w=2,4,6$, we consider finite sequences of non-negative integers $(k_3,k_4,\ldots,k_s)$
satisfying 
\begin{equation}
\sum_{j=3}^{s} k_j = m \quad\text{and}\quad 
\sum_{j=3}^{s} k_j (j-2) = w  \label{A:eqn}  \ .  
\end{equation}

\begin{kmProp}[cf.\cite{Mik:Nak}] \label{exam::a}
When  $n\ge 2$ and  weight=2, 4 or 6, 
{\color{black} the non-trivial} sub-cochain {\color{black} complexes} are as follows: 
\begin{align*} 
\cgf{1}{2n}{1}{2} & = \frakS{4} ,  
\qquad 
\cgf{2}{2n}{1}{2} 
 =  \Lambda^2 \frakS{3}  \\
 \\
\cgf{1}{2n}{1}{4}  & = \frakS{6},  
\qquad 
\cgf{2}{2n}{1}{4} 
 =  
\left(  \frakS{3} \wedge \frakS{5}  \right)
\oplus \Lambda^2 \frakS{4}  
\cong 
\left(  \frakS{3} \otimes \frakS{5}  \right)
\oplus \Lambda^2 \frakS{4},  
\\   
\cgf{3}{2n}{1}{4} 
&
 =  \Lambda^{2} \frakS{3} \wedge \frakS{4} 
 \cong  \Lambda^{2} \frakS{3} \otimes    \frakS{4}, \qquad
{\color{black}\cgf{4}{2n}{1}{4} =  \Lambda^4 \frakS{3}}
\intertext{In the above, we {\color{black} identify} the exterior product $\frakS{3} \wedge \frakS{5}$
with the tensor product $\frakS{3} \otimes  \frakS{5}$ as vector spaces, 
and we often use this identification without comments.
}
\cgf{1}{2n}{1}{6} & = \frakS{8} ,  
\qquad 
\cgf{2}{2n}{1}{6} 
 =  
\left(  \frakS{3} \otimes \frakS{7}  \right)
\oplus \left(  \frakS{4} \otimes \frakS{6}  \right)
\oplus \Lambda^2 \frakS{5}  
\\ 
\cgf{3}{2n}{1}{6} 
& = \left( \Lambda^2 \frakS{3} \otimes \frakS{6}  \right)
\oplus \left(  \frakS{3} \otimes \frakS{4} \otimes \frakS{5}  \right)
\oplus  \Lambda^{3} \frakS{4}  
\\ 
\cgf{4}{2n}{1}{6} 
& = \left( \Lambda^3 \frakS{3} \otimes \frakS{5}  \right)
\oplus \left( \Lambda^{2}  \frakS{3} \otimes \Lambda^{2} \frakS{4} \right)
\\
\cgf{5}{2n}{1}{6} 
& =  \Lambda^4 \frakS{3} \otimes \frakS{4} ,  
\qquad 
\cgf{6}{2n}{1}{6} 
 = \Lambda^{6} \frakS{3} 
\end{align*}  
\end{kmProp} 
The next job we have to do is to pick up all the trivial representations in 
$\displaystyle \cgf{m}{2n}{1}{w}$.   

\section{Crystal Base Theory}
Our group is $Sp(2n,\mR)$ and each   
irreducible representation space of $Sp(2n,\mR)$ 
is parameterized by a partition $\mu$ of length less
than $n$.  We denote its space by $\displaystyle V_{\mu}$. 
When $W$ is a representation space of $Sp(2n,\mR)$ and let us assume $W$ is
decomposed into irreducible subspaces as $\ds W = W_1 \oplus W_2 \oplus W_3$. 
If $W_1$ and $W_3$ are isomorphic with $\ds V_{\lambda}$, and 
$W_2$ is isomorphic with $\ds V_{\mu}$, then we may denote 
$\ds W \cong 2 V_{\lambda } \oplus V_{\mu} = 
2 V_{\lambda } +  V_{\mu}$ (often), and     
$\ds W^{\lambda} = W_1 \otimes W_3 \cong 2 V_{\lambda }$ and 
$\ds W^{\mu}= W_2 \cong V_{\mu} $ 

\subsection{Crystal base} 
According to \cite{kashi:french} and \cite{nakashima}, we have the crystal basis  
for $\displaystyle V_{\mu}$, 
say $\cb_{\mu}$  
consists of the semistandard {\color{\red}C}-tableaux on the
shape $\mu$, where the numbers from $1,\ldots, n, \ool{n},\ldots, \ool{1}$ are
printed, where  $\displaystyle \ool{k} := 2n+1-k$ for $ 1\le k \le n$.

For a given SS tableau, we pick a column, say J. Then J is a sequence of
increasing integers of $1,2,\ldots, n, n+1=\ool{n},\ldots, 2n = \ool{1}$.  
By $\pos{\text{J}}{a}$, we mean the position
of the element $a$ in J.  For example, if  
J = \begin{tabular}[c]{|c|}\hline 2 \\\hline   5 \\\hline 6
\\\hline\end{tabular}\ ,  
        $\pos{\text{J}}{2} = 1$, 
        $\pos{\text{J}}{5} = 2$ and  $\pos{\text{J}}{6} = 3$.   

        \newcommand{\TJ}{\text{J}}
        \newcommand{\TL}{\text{L}}
        \newcommand{\TR}{\text{R}}

SS {\color{\red}C}-tableau is a SS tableau with 2 more conditions:
\begin{itemize}
        \item[C-1] is a requirement for each column $\TJ$ of $\mu$ such that 
if 
$k, \ool{k}$ belong to $\TJ$, then 
$$\displaystyle \pos{\TJ}{k}+ (|\TJ|+1 - \pos{\TJ}{\ool{k} }) \le k\: . $$  
\item[C-2] is a requirement for each successive two columns 
        $\TL, \TR$ of $\mu$ (so that $|\TL| \ge |\TR|$) satisfy $\TL(i) \le  \TR(i)$ for
        $i=1..|\TR|$ (which is one of conditions of SS Tableau) and 
        for each $1\le a\le b\le n$,     
        if $(a,b)$-configuration then $a,b$ must satisfy  $
        \pos{}{b}-\pos{}{a} +
        \pos{}{\ool{a}} - \pos{}{\ool{b}} < b-a$.     
\end{itemize}
We assume $1 \le a \le b \le n$ (and so $1 \le a \le b < \ool{b} \le \ool{a}
\le 2n$) and assume that $a\in \TL$, $\ool{a}\in \TR$ and {$b$, $\ool{b}$} belongs
to $\TL$ or $\TR$.

$(a,b)$-configuration is when  
$$ \pos{\TL}{a} \le \pos{\TR}{b} < \pos{\TR}{\ool{b}} \le \pos{\TR}{\ool{a}} \quad\text{or}\quad 
\pos{\TL}{a} \le \pos{\TL}{b} < \pos{\TL}{\ool{b}} \le \pos{\TR}{\ool{a}} $$ 
is satisfied.      


\begin{exam} When $n=2$ and $\mu=(1,1)$, the $\cb_{\mu}$ consists of 
        \begin{center}
\begin{tabular}[c]{|c|}\hline 1\\\hline 2\\\hline\end{tabular}\hspace{5mm}%
\begin{tabular}[c]{|c|}\hline 1\\\hline 3\\\hline\end{tabular}\hspace{5mm}%
\begin{tabular}[c]{|c|}\hline 2\\\hline 3\\\hline\end{tabular}\hspace{5mm}%
\begin{tabular}[c]{|c|}\hline 2\\\hline 4\\\hline\end{tabular}\hspace{5mm}%
\begin{tabular}[c]{|c|}\hline 3\\\hline 4\\\hline\end{tabular}
\end{center}
We remark that 
\begin{tabular}[c]{|c|}\hline 1\\\hline 4=\ol{1} \\\hline\end{tabular} is a SSTab, 
but not a member of $\cb_{\mu}$ because of the condition C-1.  

The number of 
$\cb_{\mu}$ is equal to $\dim V_{\mu}$ gotten by Dimension formula coming
from Weyl's character formula.    
\end{exam}

\begin{exam} When $n=2$ and $\nu=(2,2)$, the $\cb_{\nu}$ consists of 
\begin{center} 
\begin{tabular}[c]{|c|c|}\hline 1&  1\\\hline 2&2 \\\hline\end{tabular}\hspace{3mm}%
\begin{tabular}[c]{|c|c|}\hline 1&  1\\\hline 2&3\\\hline\end{tabular}\hspace{3mm}%
\begin{tabular}[c]{|c|c|}\hline 1&  2\\\hline 2&3\\\hline\end{tabular}\hspace{3mm}%
\begin{tabular}[c]{|c|c|}\hline 1&  2\\\hline 2&4\\\hline\end{tabular}\hspace{3mm}%
\begin{tabular}[c]{|c|c|}\hline 1&  3\\\hline 2&4\\\hline\end{tabular}\hspace{3mm}%
\begin{tabular}[c]{|c|c|}\hline 1&  1\\\hline 3&3\\\hline\end{tabular}\hspace{3mm}%
\begin{tabular}[c]{|c|c|}\hline 1&  2\\\hline 3&3\\\hline\end{tabular}\hspace{3mm}%
        \\[4mm]
\begin{tabular}[c]{|c|c|}\hline 1&  2\\\hline 3&4\\\hline\end{tabular}\hspace{3mm}%
\begin{tabular}[c]{|c|c|}\hline 1&  3\\\hline 3&4\\\hline\end{tabular}\hspace{3mm}%
\begin{tabular}[c]{|c|c|}\hline 2&  2\\\hline 3&4\\\hline\end{tabular}\hspace{3mm}%
\begin{tabular}[c]{|c|c|}\hline 2&  3\\\hline 3&4\\\hline\end{tabular}\hspace{3mm}%
\begin{tabular}[c]{|c|c|}\hline 2&  2\\\hline 4&4\\\hline\end{tabular}\hspace{3mm}%
\begin{tabular}[c]{|c|c|}\hline 2&  3\\\hline 4&4\\\hline\end{tabular}\hspace{3mm}%
\begin{tabular}[c]{|c|c|}\hline 3&  3\\\hline 4&4\\\hline\end{tabular} 
\end{center} 

From the condition (C-1), each column must be one of Example 3.1.  Accounting
SS Tab property, then candidates should be $\displaystyle {5-1+2 \choose 2} =
15$, but 
\begin{tabular}[c]{|c|c|}\hline 2&  2\\\hline 3&3\\\hline\end{tabular}\ 
has $(2,2)$-configuration property but does not satisfy the condition (C-2),
and is not a member of $\cb_{\nu}$.  
\end{exam}

\subsection{Decomposition of tensor product by crystal base} 

Take another irreducible representation $\displaystyle V_{\lambda}$. Then the
irreducible decomposition of the tensor
product  $\displaystyle V_{\lambda} \otimes V_{\mu}$ is given by
combinatorially from $\lambda$ and the crystal base of $\mu$ as follows:  

For each tableau $T$ of $\cb_{\mu}$, 
\begin{itemize}
\item[(1)] Stand at the top of the right most column of $T$.  

\textbf{action:} 
from the top to downward, let the
number on the cell apply to the Young diagram by the next rule: 

\begin{quote}
    If the number $k$ on the cell of $T$ satisfies  
    $k \le n$, then  
    the new diagram is made by adding one cell to the $k$-th row 
    of the old Young diagram. If the new diagram is not Young diagram, stop
    this process.  

    If the number $k$ satisfies $k > n$, then 
    the new diagram is made by deleting one cell from
    the $(2n+1-k)$-th row of the old Young diagram. 
    If the new diagram is not Young diagram, stop this process.  
\end{quote}
\item[(2)] If we still have a Young diagram, then move the one left column and do 
\textbf{the action} above, and so on.   

\item[(3)] 
		At the left most column, 
		if we could get finally a Young diagram, say $\nu$, 
		after the all actions, 
		then $\displaystyle
V_{\nu}$ is one factor of $V_{\lambda}\otimes V_{\mu}$.  
\end{itemize}

\subsubsection{A concrete example of the action of 
crystal base to a Young diagram}
When $n=3$, the crystal base of 
$\mu = (2,1)$ 
=
\begin{tabular}[c]{*{2}{|c}|}
        \hline
          &   \\ \hline  
          \\
        \cline{1-1}
\end{tabular}
is known as 
\begin{center} 
\begin{tabular}[c]{*{2}{|c}|}
        \hline 1 & 1  \\ \hline  
               2 \\ \cline{1-1} \end{tabular}, 
\begin{tabular}[c]{*{2}{|c}|}
        \hline 1 & 1  \\ \hline  
               3 \\ \cline{1-1} \end{tabular}, 
\begin{tabular}[c]{*{2}{|c}|}
        \hline 1 & 1  \\ \hline  
               \ol{3} \\ \cline{1-1} \end{tabular}, 
\begin{tabular}[c]{*{2}{|c}|}
        \hline 1 & 1  \\ \hline  
               \ol{2} \\ \cline{1-1} \end{tabular}, 
\begin{tabular}[c]{*{2}{|c}|}
        \hline 1 & 1  \\ \hline  
               \ol{1} \\ \cline{1-1} \end{tabular}, 
               \ldots 
\begin{tabular}[c]{*{2}{|c}|}
        \hline 3 & \ol{1}  \\ \hline  
               \ol{2} \\ \cline{1-1} \end{tabular}, 
\begin{tabular}[c]{*{2}{|c}|}
        \hline 3 & \ol{1}  \\ \hline  
               \ol{1} \\ \cline{1-1} \end{tabular}, 
\begin{tabular}[c]{*{2}{|c}|}
        \hline \ol{3} & \ol{1}  \\ \hline  
               \ol{2} \\ \cline{1-1} \end{tabular}, 
\begin{tabular}[c]{*{2}{|c}|}
        \hline \ol{3} & \ol{1}  \\ \hline  
               \ol{1} \\ \cline{1-1} \end{tabular}, 
\begin{tabular}[c]{*{2}{|c}|}
        \hline \ol{2} & \ol{1}  \\ \hline  
               \ol{1} \\ \cline{1-1} \end{tabular} 
       \end{center} 

Apply \begin{tabular}[c]{*{2}{|c}|}
        \hline 1 & 1  \\ \hline  
               2 \\ \cline{1-1} \end{tabular}\quad to 
$\lambda=(1,1) = 
\begin{tabular}[c]{*{1}{|c}|} \hline 
          \\ \hline  
          \\ \cline{1-1} \end{tabular}\quad 
$ then  we have 
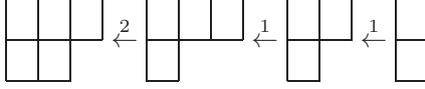
\begin{center} 
\begin{tabular}[c]{*{3}{|c}|} \hline 
        & &  \\ \hline  
        &   \\ \cline{1-2} \end{tabular} 
$\displaystyle \mathop{\leftarrow}^{2} $ 
\begin{tabular}[c]{*{3}{|c}|} \hline 
        & &  \\ \hline  
             \\ \cline{1-1} \end{tabular}
$\displaystyle \mathop{\leftarrow}^{1} $ 
\begin{tabular}[c]{*{2}{|c}|} \hline 
         &   \\ \hline  
             \\ \cline{1-1} \end{tabular} 
$\displaystyle \mathop{\leftarrow}^{1} $ 
\begin{tabular}[c]{*{1}{|c}|} \hline 
          \\ \hline  
          \\ \cline{1-1} \end{tabular} 
               \end{center}

Apply \begin{tabular}[c]{*{2}{|c}|}
        \hline 1 & 1  \\ \hline  
\ol{3}  \\ \cline{1-1} \end{tabular}\quad to 
$\lambda=(1,1) = 
\begin{tabular}[c]{*{1}{|c}|} \hline 
          \\ \hline  
          \\ \cline{1-1} \end{tabular}\quad 
$ then by the notation of partitions, we have  
\begin{center} 
   Not Young diagram     $\displaystyle (3,1,-1) \mathop{\leftarrow}^{\ool{3}}  
        (3,1) \mathop{\leftarrow}^{1}   (2,1) \mathop{\leftarrow}^{1}  (1,1) $
\end{center}

Apply \begin{tabular}[c]{*{2}{|c}|}
        \hline 1 & 3  \\ \hline  
\ol{3}  \\ \cline{1-1} \end{tabular}\quad to 
$\lambda=(1,1) = 
\begin{tabular}[c]{*{1}{|c}|} \hline 
          \\ \hline  
          \\ \cline{1-1} \end{tabular}\quad 
$ then by the notation of partitions, we have  
\begin{center} 
   $\displaystyle (2,1,0) \mathop{\leftarrow}^{\ool{3}}  
 (2,1,1) \mathop{\leftarrow}^{1} (1,1,1) \mathop{\leftarrow}^{3} (1,1,0) $
\end{center}

Apply \begin{tabular}[c]{*{2}{|c}|}
		\hline 1 & \ol{2}  \\ \hline  
2  \\ \cline{1-1} \end{tabular}\quad to 
$\lambda=(1,1) = 
\begin{tabular}[c]{*{1}{|c}|} \hline 
          \\ \hline  
          \\ \cline{1-1} \end{tabular}\quad 
$ then by the notation of partitions, we have  
\begin{center} 
   $\displaystyle (2,1,0) \mathop{\leftarrow}^{2}  
   (2,0,0)\mathop{\leftarrow}^{1} (1,0,0)\mathop{\leftarrow}^{\ool{2}} (1,1,0)$
\end{center}

Apply \begin{tabular}[c]{*{2}{|c}|}
        \hline \ol{2} & \ol{1}  \\ \hline  
\ol{1}  \\ \cline{1-1} \end{tabular}\quad to 
$\mu=(2,1) = 
\begin{tabular}[c]{*{2}{|c}|} \hline 
        & \\ \hline  
          \\ \cline{1-1} \end{tabular}\quad  
$ then by the notation of partitions, we have 
\begin{center} 
 trivial $\displaystyle  (0,0) \mathop{\leftarrow}^{\ool{1}}  
   (1,0) \mathop{\leftarrow}^{\ool{2}}  
   (1,1) \mathop{\leftarrow}^{\ool{1}}  
        (2,1)$.   
\end{center}

\subsection{Pick up trivial representations from $\displaystyle
V_{\lambda}\otimes V_{\mu}$}
It is well-known that $triv \otimes V_{\mu} \cong V_{\mu}$ for each $\mu$. If
we look at this result by the crystal base theory, we see that 
there is the unique element $J$ of the crystal base 
of $\displaystyle V_{\mu}$ and satisfies 
$J\cdot triv = \mu$,
where $\cdot$ means the action of $J$ to each partition.  

\begin{kmProp}\label{prop:triv:first}
The SS {\color{\red}C}-tableau $J$ of the shape $\mu$, satisfying 
		$\displaystyle J \cdot triv = \mu$ is defined by 
printing number $k$ for the cells of the $k$-row of the shape of $\mu$.  
We denote this SS {\color{\red}C}-tableau by $T_{\mu}$.   

If a SS {\color{\red}C}-tableau $J$ of the shape $\mu$ satisfies $\displaystyle
J\cdot triv = \lambda$, then $\lambda=\mu$ and $\displaystyle J = T_{\mu}$.  If
$\mu \ne triv$, then at each step of $J\cdot triv$, the Young diagram is not
trivial.  

\end{kmProp}
\textbf{Proof:} The first operation acting to the trivial Young diagram is $1$. Thus, the right most and the top entry of $J$ is $1$.  By SS Tab properties, 
the first row of $J$ is printed by only $1$.  
Now by (C1) condition, $J$ does have no \ol{1}.  
Let $x$ be the entry of left most of the second row of $J$.  Then the first part of the sequence of operations is
$\displaystyle x,1,\underbrace{1,\ldots,1}_{(\mu_1 - \mu_2)\text{-times}}$
with $ x \ne \ool{1}$. Keeping the Young diagram property implies $x = 2$. We claim that 
$J$ does have no \ol{2}.  The reason is: if $\ool{2}\in J$, take a column $L$ of $J$ with $\ool{2}\in L$. Then $\displaystyle \pos{L}{2}=2$ and 
$\displaystyle \pos{L}{\ol{2}} = |L|$ and conflicts with (C1) condition. Repeating the same discussion, we concludes our proof. \kmqed

We denote the diagram of height $h$ of one cell by 
$s_h$. Any Young diagram of length less than $n$ is expressed as 
$$\mu = s_n{}^{p_n}  s_{n-1}{}^{p_{n-1}} \cdots  s_2{}^{p_2}  s_{1}{}^{p_1} 
$$
where $p_i \ge 0$ ($i=1..n$).     
\begin{align} \label{pat:tower}
		\mu_1 =& p_n + p_{n-1} + \cdots + p_2 + p_1 &&& p_n &= \mu_n\notag\\
		\mu_2 =& p_n + p_{n-1} + \cdots + p_2 &&& p_{n-1} &= \mu_{n-1}-\mu_n\notag\\
		\vdots &&&& \vdots & \\
		\mu_n =& p_n &&& p_1 &= \mu_1 - \mu_2\notag 
\end{align}
holds good. Conversely, for a given partition $\mu=(\mu_1\ge \mu_2\ge \cdots
\ge \mu_{n} \ge 0)$, we define the sequence of $\displaystyle p_i$  
($i=1..n$) the above relations.  Thus, sometimes we regard of Young diagram as
the (horizontal!?) collection of columns.  

Let $\displaystyle t_{j} = 
\begin{tabular}[c]{|c|}
        \hline
        1 \\ 2 \\
        \vdots  \\ j\\\hline
\end{tabular}$ be a SS {\color{\red}C}-tableau on $\ds s_j$. 
Then $\ds T_{\mu} = 
 t_n{}^{p_n}  t_{n-1}{}^{p_{n-1}} \cdots  t_2{}^{p_2}  t_{1}{}^{p_1}$.

\begin{defn}
		For a partition or Young diagram  $\mu$ 
		of length less than $n$, define the 
		SS {\color{\red}C}-tableau $\displaystyle T_{\mu}$ by 
		printing number $k$ for the cells of the $k$-row of 
		$\displaystyle T_{\mu}$. In other words, 
$\ds T_{\mu} = 
 t_n{}^{p_n}  t_{n-1}{}^{p_{n-1}} \cdots  t_2{}^{p_2}  t_{1}{}^{p_1}$, 
 where $p_j$ ($j=1..n$) are defined as (\ref{pat:tower}).     
\end{defn}

\begin{exam}
        When $n=3$ and $\mu= (4,2,1)$, 
        $T = \begin{tabular}{*{4}{|c}|} \hline
                1 & 1 & 1 & 1 \\\hline 
                2 & 2  \\\cline{1-2} 
                3   \\\cline{1-1} 
        \end{tabular}$, then 
$ T\cdot triv $ becomes $\mu$ as $$ 
(4,2,1) \mathop{\leftarrow}^{3}
(4,2,0) \mathop{\leftarrow}^{2}
(4,1,0) \mathop{\leftarrow}^{1}
(3,1,0) \mathop{\leftarrow}^{2}
(3,0,0) \mathop{\leftarrow}^{1}
(2,0,0) \mathop{\leftarrow}^{1}
(1,0,0) \mathop{\leftarrow}^{1}
(0,0,0) \; . $$
We look the above flow reversely, i.e., 
$$ 
(4,2,1) \mathop{\rightarrow}^{\ool{3}}
(4,2,0) \mathop{\rightarrow}^{\ool{2}}
(4,1,0) \mathop{\rightarrow}^{\ool{1}}
(3,1,0) \mathop{\rightarrow}^{\ool{2}}
(3,0,0) \mathop{\rightarrow}^{\ool{1}}
(2,0,0) \mathop{\rightarrow}^{\ool{1}}
(1,0,0) \mathop{\rightarrow}^{\ool{1}}
(0,0,0)\; , $$
the action does not come from a SS {\color{\red}C}-tableau. 
Instead of this, we apply bar-operation for each element of $T$ 
and get a new tableau 
        $T' = \begin{tabular}{*{4}{|c}|} \hline
                \ol{1} & \ol{1} & \ol{1} & \ol{1} \\\hline 
                \ol{2} & \ol{2}  \\\cline{1-2} 
                \ol{3}   \\\cline{1-1} 
        \end{tabular}$\ . Although
$T'$ is not SS Tab, but reversing the elements of each column of $T'$, 
we get a SS {\color{\red}C}-tableau 
$\widehat{T} = \begin{tabular}{*{4}{|c}|} \hline
                \ol{3} & \ol{2} & \ol{1} & \ol{1} \\\hline 
                \ol{2} & \ol{1}  \\\cline{1-2} 
                \ol{1}   \\\cline{1-1} 
        \end{tabular}$\ . 
        Acting $\displaystyle \widehat{T}$ on the partition $\mu=(4,2,1)$, we get
        the trivial representation as below:  
$$ 
(0,0,0) \mathop{\leftarrow}^{\ool{1}}
(1,0,0) \mathop{\leftarrow}^{\ool{2}}
(1,1,0) \mathop{\leftarrow}^{\ool{3}}
(1,1,1) \mathop{\leftarrow}^{\ool{1}}
(2,1,1) \mathop{\leftarrow}^{\ool{2}}
(2,2,1) \mathop{\leftarrow}^{\ool{1}}
(3,2,1) \mathop{\leftarrow}^{\ool{1}}
(4,2,1)$$ 
\end{exam}
\begin{defn}
We put $\ds  \hat{t}_{j} = 
\begin{tabular}[c]{|c|}
        \hline
        \ol{j} \\ \vdots \\
        \ol{2} \\ \ol{1}\\\hline
\end{tabular}$ for $j=1..n$.  Let $\mu$ be a partition of length less than or equal to $n$.  Define a tableau by 
$\ds \widehat{T}_{\mu} = 
 \hat{t}_n{}^{p_n}  \hat{t}_{n-1}{}^{p_{n-1}} \cdots  \hat{t}_2{}^{p_2}
 \hat{t}_{1}{}^{p_1}$,  where $p_j$ ($j=1..n$) are defined as in (\ref{pat:tower}).     
 \end{defn}
 It is easy to see that $\ds \widehat{T}_{\mu}$ is a SS {\color{\red}C}-tableau 
 on the shape $\mu$ and satisfies 
 $\ds \widehat{T}_{\mu} \cdot \mu = triv$ because of  
\begin{align*}
		\widehat{T}_{\mu} \cdot \mu =& 
		( \hat{t}_n{}^{p_n}  \hat{t}_{n-1}{}^{p_{n-1}} \cdots  \hat{t}_2{}^{p_2}  \hat{t}_{1}{}^{p_1}) \cdot ( 
		s_n{}^{p_n}  s_{n-1}{}^{p_{n-1}} \cdots  s_2{}^{p_2}  s_{1}{}^{p_1}) \\
		=&
		( \hat{t}_n{}^{p_n}  \hat{t}_{n-1}{}^{p_{n-1}} \cdots \hat{t}_2{}^{p_2})\cdot   \hat{t}_{1}{}^{p_1}) \cdot ( 
		s_n{}^{p_n}  s_{n-1}{}^{p_{n-1}} \cdots  s_2{}^{p_2}  s_{1}{}^{p_1}) \\
		=&
		( \hat{t}_n{}^{p_n}  \hat{t}_{n-1}{}^{p_{n-1}} \cdots \hat{t}_2{}^{p_2})\cdot   ( 
		s_n{}^{p_n}  s_{n-1}{}^{p_{n-1}} \cdots  s_2{}^{p_2}  ) \\
		\vdots \\
		=& triv	
\end{align*}

As reverse situations of Proposition \ref{prop:triv:first}, we have  
\begin{Lemma}\label{key::lemma} Let   
$J$ be a member of the crystal base of 
the irreducible representation 
$V_{\mu}$ of $Sp(2n,\mR)$, i.e. $J$ is a SS {\color{\red}C}-tableau on the
shape $\mu$.  If  $J\cdot \lambda = triv$, then $\lambda = \mu$ and $\ds J=
\widehat{T}_{\mu}$.  
\end{Lemma}
\textbf{Proof:} 
Suppose $\ds J\cdot \lambda = triv$. Then the last operation of $J$ is \ol{1}.
This means $J$ has \ol{1} on the bottom of the left most column, say $L$.  Then
it turns out $J$ does not have the number 1. The reason is that if $J$ has 1,
then $L$ can be characterized as the top entry is 1 and the bottom entry is
\ol{1}. This conflicts with the condition (C-1).  

$J$ having no 1 means that on the way of successive actions $\ds J\cdot
\lambda$, the trivial Young diagram does not appear until the final stage.    

Let us imagine the one more past shape of the Young diagram.  
The possibilities are two, 
$\ds triv \mathop{\leftarrow}^{\ool{1}}  \  
\begin{tabular}{|c|} \hline \\\hline\end{tabular}
		\ds \mathop{\leftarrow}^{\ool{1}}  
		\begin{tabular}{*{2}{|c}|}\hline & \\\hline \end{tabular}\ $ 
or 
$\ds triv \mathop{\leftarrow}^{\ool{1}} \  
\begin{tabular}{|c|} \hline \\\hline\end{tabular}
		 \mathop{\leftarrow}^{\ool{2}}  
\begin{tabular}{|c|} \hline  \\\hline \\\hline\end{tabular}\ $.    
In the first case,  $J$ becomes 
\begin{tabular}{*{3}{|c}|}\hline \ol{1} & $\cdots$ & \ol{1}\\\hline\end{tabular}\ .  

In the second case,  $L$ is in form of 
\begin{tabular}{|c|} $\vdots$ \\\hline \ol{2}\\\hline
\ol{1}\\\hline \end{tabular}\ , 
and we claim that the number 2 does not live in  $J$ because of (C-1)
condition.  

Repeating the same arguments, we may assume that there is a maximal number 
$k$ satisfying $ \ool{j}\in L$ but  $\ds j\not\in J$ for $j=1..k$.

Again, imagine the Young diagram one step before.  The three possibilities: 
$\ds [1^k] \leftarrow [1^{k-1}]$, 
$\ds [1^k] \leftarrow [2, 1^{k-1}]$, or 
$\ds [1^k] \leftarrow [1^{k+1}]$.  The first case is $y=k$ (impossible), 
the second case is $y=\ool{1}$ and the third case is $y = \ool{k+1}$ (impossible). From those, we conclude $\ds L = \hat{t}_{k}$ and if 
$J$ has the second column from the left, say $M$,  
the bottom entry of $M$ 
is \ol{1} and is just 
the successive series of \ol{1}, \ol{2} \ldots upward. We continue the discussion step by step and reach to the right most column. 
Thus, 
it turns out $J = \widehat{T}_{\mu}$.

We write down $\lambda$ as
$\ds \lambda = s_n{}^{q_n} \cdots s_2{}^{q_2} s_1{}^{q_1}$, and 
consider the contributions of $\ds\widehat{T}_{\mu}$ to $\lambda$, then 
\begin{align*}
		\widehat{T}_{\mu} \cdot \lambda =& 
(\hat{t}_n{}^{p_n}\hat{t}_{n-1}{}^{p_{n-1}}\cdots\hat{t}_2{}^{p_2}\hat{t}_{1}{}^{p_1})\cdot 
(s_n{}^{q_n}s_{n-1}{}^{q_{n-1}}\cdots s_2{}^{q_2} s_{1}{}^{q_1}) \\
		=& 
(\hat{t}_n{}^{p_n}\hat{t}_{n-1}{}^{p_{n-1}}\cdots\hat{t}_2{}^{p_2})\cdot  
(s_n{}^{q_n}s_{n-1}{}^{q_{n-1}}\cdots s_2{}^{q_2} s_{1}{}^{q_1 -p_1}) 
\quad \text{if } q_1 - p_1 \ge 0
\\ =& 
(\hat{t}_n{}^{p_n}\hat{t}_{n-1}{}^{p_{n-1}}\cdots\hat{t}_3{}^{p_3})\cdot 
(s_n{}^{q_n}s_{n-1}{}^{q_{n-1}}\cdots s_2{}^{q_2-p_2 } 
s_1{}^{ q_1 - p_1 }) 
\quad \text{if } q_1 - p_1 \ge 0, q_2 - p_2 \ge 0 
\\ \vdots & \\
=& 
(s_n{}^{q_n - p_n}\cdots s_3{}^{q_3 -p_3} 
s_2{}^{q_2 -p_2} 
s_1{}^{q_1 -p_1} 
) 
\quad \text{if } q_j - p_j \ge 0 \text{ for } j = 1..n
\end{align*}
and we conclude 
$ \mu_j = \lambda_j$ ( $j=1..n$), namely $\ds \lambda = \mu$.    
\kmqed

\begin{kmProp}

Let $W$ and $Z$ be $Sp(2n,\mR)$-representation space and assume that they have
irreducible decompositions like 
$\ds 
W = \mathop{\oplus}_{\lambda\in I} \xi_{W}^{\lambda} V_{\lambda} \quad  
 (\xi_{W}^{\lambda} \in \mathbb{N})$ and   
 $\ds 
Z = \mathop{\oplus}_{\mu\in J} \xi_{Z}^{\mu} V_{\mu} 
\quad ( \xi_{Z}^{\mu} \in \mathbb{N}) $.   
Then the tensor product 
$W\otimes Z$ of $W$ and $Z$ has the trivial representation of 
the multiplicity 
$\ds 
\sum_{\lambda\in I\cap J} 
\xi_{W}^{\lambda} 
\xi_{Z}^{\lambda} $, i.e., 
$$  \left( W \otimes Z\right)^{triv} \cong  
\left(\sum_{\lambda\in I\cap J} 
\xi_{W}^{\lambda} 
\xi_{Z}^{\lambda}\right) \  triv 
$$ 
\end{kmProp}

\bigskip

\begin{kmCor}
Let $W$  be a $Sp(2n,\mR)$ representation and assume that $W$ has an irreducible
decomposition as 
$\displaystyle 
W = \mathop{\oplus}_{\lambda\in I} \xi_{W}^{\lambda} V_{\lambda} \quad  
 (\xi_{W}^{\lambda} \in \mathbb{N})$.    
Then we have 
$$ 
\dim\; \left( W \otimes \frakS{h} \right)^{triv} = 
\xi_{W}^{[h,0^{n-1}] } $$ 
\end{kmCor}

\section{The relative cochain complexes}
When  $n\ge 2$ and  weight=2, 4 or 6, we have 
{\color{black} the non-trivial} sub-cochain {\color{black} complexes} 
$\ds \cgf{\bullet}{2n}{1}{w}$   
in Proposition \ref{exam::a}.   
We apply elementary irreducible decomposition rules to those spaces
$\ds \cgf{\bullet}{2n}{1}{w}$ and have the   
$\ds \CGF{\bullet}{2n}{0}{w}$  as follows.  
\begin{align*} 
\CGF{1}{2n}{0}{2} & = 0 ,  
        \hspace{35mm}
\CGF{2}{2n}{0}{2} = \left(\Lambda^2 \frakS{3}\right)^{triv}  \\
 \\
\CGF{1}{2n}{0}{4}  & = 0 ,  
 \hspace{35mm}
\CGF{2}{2n}{0}{4} 
\cong \left(\Lambda^2 \frakS{4}\right)^{triv},  
\\   
\CGF{3}{2n}{0}{4} 
&
= \left(( \Lambda^{2} \frakS{3} ) \otimes \frakS{4}\right)^{triv} , \hspace{5mm}
{\color{black}\CGF{4}{2n}{0}{4} =  \left(\Lambda^4 \frakS{3}\right)^{triv}}
\\
\\
\CGF{1}{2n}{0}{6} & = 0  ,  
\hspace{35mm}
\CGF{2}{2n}{0}{6} 
 =  \left(\Lambda^2 \frakS{5}\right)^{triv}\; ,  
\\ 
\CGF{3}{2n}{0}{6} 
& = \left( (\Lambda^2 \frakS{3}  ) \otimes \frakS{6}\right)^{triv} 
\oplus  \left(\frakS{3} \otimes \frakS{4} \otimes \frakS{5}\right)^{triv}
\oplus \left( \Lambda^{3} \frakS{4}\right)^{triv}\; ,   
\\ 
\CGF{4}{2n}{0}{6} 
& = \left( (\Lambda^3 \frakS{3} ) \otimes \frakS{5}\right)^{triv} 
\oplus \left(\Lambda^{2}\frakS{3}\otimes\Lambda^{2}\frakS{4}\right)^{triv}\; ,
\\
\CGF{5}{2n}{0}{6} 
& =  \left( (\Lambda^4 \frakS{3}) \otimes \frakS{4} \right) ^{triv}   ,  
\hspace{5mm}
\CGF{6}{2n}{0}{6} 
= \left(\Lambda^{6} \frakS{3}\right)^{triv} 
\end{align*}  

\subsection{Splitting $\ds\Lambda^p\frakS{q}$ when $n=3$}
As we see above lists, we need some irreducible factors or complete
irreducible decomposition of $\ds \Lambda^{p}\frakS{q}$ 
case by case,   
which is 
an $Sp(6,\mR)$-invariant subspace of $\ds
\mathop{\otimes}^{p} \frakS{q}$.       

We recall here the fundamental fact which our discussion is completely based
on.  We state the fact for $n$=3 even though it is valid for a general $n$. 

Let $T$ be the subgroup of $Sp(6,\mR)$ consisting of diagonal matrices, 
$U$ be the subgroup of $Sp(6,\mR)$ consisting of upper triangle matrices whose
diagonal entries are 1, and  
$U^{-} $ be the subgroup of $Sp(6,\mR)$ consisting of lower triangle matrices whose
diagonal entries are 1.  $U$ is called the maximal unipotent subgroup of
$Sp(6,\mR)$.    

Let $W$ be a representation of $Sp(6,\mR)$. $\mathbf{x} \in W $ is called a
maximal vector if it satisfies 
\begin{equation}\label{max:vec} \mathbf{x}\ne \mathbf{0},\quad  g\cdot \mathbf{x}
        = \mathbf{x}\quad (\forall g\in U), \quad 
t\cdot \mathbf{x}
= \lambda(t) \mathbf{x} \quad  
(\forall t\in T, \ \exists \text{ weight } \lambda)\ .
\end{equation}
The maximal vector is unique up to scalar multiple.  
\medskip

\textbf{Fact:} Let $W$ be a representation of $Sp(6,\mR)$. Then the next hold.
\begin{itemize}
\item[(1)] Let $\mathbf{x}$ a maximal vector. Then the $Sp(6,\mR)$-invariant subspace
of $W$ generated by $\mathbf{x}$ is an irreducible subspace. The irreducible
space is also generated by $\ds \{ (U^{-})^{k} \cdot \mathbf{x} \mid
k=0,1,2,\ldots\}$.  
\item[(2)] $W$ can be decomposed into irreducible subspaces which are one-to-one
        corresponding to the set of maximal vectors.
\end{itemize}

\bigskip

To find maximal vectors, we consider the infinitesimal version of the second
condition of (\ref{max:vec}), namely,  
$\ds \xi \cdot \mathbf{x} = \mathbf{0}$ for 
($\forall \xi\in \frak{u}$), where 
$\frak{u}$ is the Lie algebra of $U$ and is 9-dimensional.  

The infinitesimal action $\xi\in\frak{u} $ 
behaves as an ordinary derivation of degree 0 and for each 1-cochain $\sigma$
and for each polynomial $f$, we have 
$$ \langle \xi \cdot \sigma , f \rangle = - \langle \sigma , 
\xi \cdot f\rangle =
- \langle \sigma, \Pkt{ \hat{J}(\xi)}{f}\rangle  
$$
where $\ds J$ is the momentum mapping of the natural symplectic action and 
$\ds \Pkt{}{}$ is the standard Poisson bracket on  $\ds\mR^6$.    
In general, we have to deal with a big simultaneous linear homogeneous
equations of 
$9\times \dim W$-equations with 
$\dim W$-variables.

Since our spaces are 
concrete and special, namely they are (sum of wedge products of) 
dual of homogeneous polynomials, 
we devise a simple but
effective treatment for 
the third condition of (\ref{max:vec}).   
By the Crystal Base theory we know the decomposition of 
$\ds \mathop{\otimes}^{p}\frakS{q}$ and this decomposition suggests all the
possibilities of maximal weights which appear in $\ds \Lambda^p\frakS{q}$. 

We explain this trick in the case of $\ds\Lambda^{2}\frakS{5}$.  If we write
each basic element of $\ds\frakS{5}$ by $\ds z_A$ where $\ds A =
(a_1,a_2,\ldots,a_6)$ with $a_j \in \mN^{+}$ and $\ds\sum_{j=1}^6 a_j = 5$.  
A basis of 
$\ds\Lambda^{2}\frakS{5}$ consists of $\ds z_A \wedge z_B$ with $A <B$ in
the lexicographic order.  We see that 
$\ds\dim\frakS{5}=252$ when $n=3$ and so $\ds\dim\Lambda^{2}\frakS{5}=
31626$. Observing the torus action of $\ds Sp(6,\mR)$ on this space, 
we see that
$$[ (a_1 - a_6) + (b_1 - b_6),  (a_2 - a_5) + (b_2 - b_5),  
(a_3 - a_4) + (b_3 - b_4)]$$  means the weight of 
$\ds z_A \wedge z_B$. If we try to find some basis of the trivial space, then 
the candidates are 
$\ds z_A \wedge z_B$ with 
$[ (a_1 - a_6) + (b_1 - b_6),  (a_2 - a_5) + (b_2 - b_5),  
(a_3 - a_4) + (b_3 - b_4)]=[0,0,0]$ and we see the dimension is 330. 
We stress that the number is drastically reduced. 
We prepare a general vector of the form $\ds v = \sum_{j=1}^{330} c_j w_j$,
where $\ds w_j $ is one of $\ds z_A \wedge z_B$ with the requirement above. 
By invariance by the maximal unipotent subgroup $U$, 
we have 2159 linear equations of $c_j$. Solving the linear equations, we
see that the freedom is 1, namely the solution space is 1-dimensional, namely 
$\ds\Lambda^2\frakS{5}$ has the trivial representation space with multiplicity
1.    
If we restart from another weight $\ds \lambda = [\lambda_1,\lambda_2,\lambda_3]$, 
we may continue the same
discussion and get the multiplicity of $\ds V_{\lambda}$ in
$\ds\Lambda^2\frakS{5}$. 
This discussion works well for general 
$\ds\Lambda^p\frakS{q}$. 

\bigskip

We show the irreducible decompositions of some $\ds
\Lambda^p \frakS{q}$ without proof.

\begin{kmProp}\label{prop:lambda:s}
When $n=3$, we have 
\begin{align} 
				\Lambda^2 \frakS{3}  \cong &
V_{[0,0,0]}+V_{[1,1,0]}+V_{[2,2,0]}+V_{[3,3,0]}+V_{[4,0,0]}+V_{[5,1,0]} 
\label{ni:san}
                 \\
				\Lambda^2 \frakS{4}  \cong &
 V_{[2,0,0]}+V_{[3,1,0]}+V_{[4,2,0]}+V_{[5,3,0]}+V_{[6,0,0]}+V_{[7,1,0]} 
 \label{ni:yon}
 \\
				\Lambda^2 \frakS{5}  \cong &
 V_{[0,0,0]}+V_{[1,1,0]}+V_{[2,2,0]}+V_{[3,3,0]}+V_{[4,0,0]}+V_{[4,4,0]}\notag\\&
 +V_{[5,1,0]}+V_{[5,5,0 ]}
                 +V_{[6,2,0]}+V_{[7,3,0]}+V_{[8,0,0]}+V_{[9,1,0]} 
                 \label{ni:go}
\\
\Lambda^3\frakS{3} =& V_{[2,1,0]}+3 V_{[3,0,0]}+2
V_{[3,1,1]}+V_{[3,2,0]}+V_{[3,2,2]}+V_{[3,3,3]}+2 V_{[4,1,0]}\notag\\&
+V_{[4,2,1]} +V_{[4,3,0]}+2 V_{[5,2,0]}+V_{[5,3,1]}+V_{[6,1,0]}
+V_{[6,3,0]}+V_{[7,0,0]}+V_{[7,1,1]} 
\label{san:san}
\\
        \Lambda^3\frakS{4} =& 
        2 V_{[2,0,0]}+V_{[2,1,1]}+V_{[2,2,2]}+3 V_{[3,1,0]}+V_{[3,2,1]}+V_{[3,3,0]}
        +V_{[3,3,2]}+V_{[4,0,0]}\notag\\& 
        +V_{[4,1,1]}+4 V_{[4,2,0]}+2 V_{[4,3,1]}+V_{[4,4,2]}
        +3 V_{[5,1,0]}+2 V_{[5,2,1]}+3 V_{[5,3,0]}+V_{[5,3,2]}\notag\\& +V_{[5,4,1]}
        +V_{[5,5,2]}
       +3 V_{[6,0,0]}+2 V_{[6,1,1]}+2
       V_{[6,2,0]}+V_{[6,2,2]}+V_{[6,3,1]}+V_{[6,3,3]}\notag\\&
       +2 V_{[6,4,0]}+2 V_{[7,1,0]}+V_{[7,2,1]}+2 V_{[7,3,0]}+V_{[7,4,1]}+V_{[7,5,0]}
       +2 V_{[8,2,0]}+V_{[8,3,1]}\notag\\&
       +V_{[9,1,0]}+V_{[9,3,0]}+V_{[10,0,0]}+V_{[10,1,1]}
       \label{san:yon}
      \\ 
        \Lambda^4\frakS{3} =&  
3 V_{[0,0,0]}+4 V_{[1,1,0]}+V_{[2,1,1]}+6 V_{[2,2,0]}+3 V_{[3,1,0]} +4
V_{[3,2,1]}+7 V_{[3,3,0]} +V_{[3,3,2]}\notag\\& +4 V_{[4,0,0]} +5 V_{[4,1,1]} 
+4 V_{ [4,2,0]}+3 V_{[4,2,2]}+4 V_{[4,3,1]}+V_{[4,3,3]}+4 V_{[4,4,0]}+7
V_{[5,1 ,0]}\notag\\&
+5 V_{[5,2,1]}+2 V_{[5,3,0]}
+3 V_{[5,3,2]}+2 V_{[5,4,1]}+2 V_{[5,5,0 ]}+2 V_{[6,0,0]}+3 V_{[6,1,1]}+5 V_{[6,2,0]}
\notag \\&
+V_{[6,2,2]}+3 V_{[6,3,1]}+V_{[6,3,3]}
+V_{[6,4,0]}+V_{[6,4,2]}+V_{[6,6,0]}+V_{[7,1,0]}+2 V_{[7,2,1]}
\notag\\&
+
3 V_{[7,3,0]} 
+V_{[7,4,1]}+V_{[8,0,0]}+V_{[8,1,1]} +V_{[8,2,0]}+V_{[8,3,1]}+V_{[9,1,0]} 
\label{yon:san}
\end{align}
\end{kmProp}
\begin{kmRemark} 
As we mentioned in \cite{Mik:Nak} when $n=2$, we know the dimensions of each
irreducible components above by Weyl's dimension formula and we compare the
both sides of the equation above in dimension, and we can check validity of
our decomposition by comparing the dimensions.    
\end{kmRemark}
\subsection{The dimensions of relative cochain complexes}
We abbreviate the space   
 $\ds \CGF{m}{6}{0}{w}$ by $\ds \CCx{m}{w}$.  
Applying Proposition \ref{prop:lambda:s}, we see that 
\begin{align*}
\dim \CCx{2}{2} &= \dim \left(\Lambda^2\frakS{3}\right)^{triv}
 \mathop{=}^{ (\ref{ni:san})} 1 \\ 
\dim \CCx{2}{4} &= \dim \left(\Lambda^2\frakS{4}\right)^{triv}
 \mathop{=}^{ (\ref{ni:yon})} 0 \\ 
 \dim \CCx{3}{4} &= \dim \left((\Lambda^2\frakS{3}) \otimes \frakS{4} \right)^{triv}
 = \dim \left( \Lambda^2 \frakS{3}\right)^{[4,0,0]} 
 \mathop{=}^{ (\ref{ni:san})} 1 \\ 
\dim \CCx{4}{4} &= \dim \left(\Lambda^4\frakS{3}\right)^{triv}
 \mathop{=}^{ (\ref{yon:san})} 3 
 \end{align*} 
Thus, for weight =2 and 4 we get the dimensions of each cochain complexes.  
\begin{center}
        \begin{tabular}{c|*{2}c}
                & $\ds\CCx{1}{2}$ & $\CCx{2}{2}$ \\\hline
                $\dim$ & 0 & 1
        \end{tabular}
\qquad and 
\qquad 
        \begin{tabular}{c|*{4}c}
                & $\ds\CCx{1}{4}$ & $\CCx{2}{4}$ & $\ds\CCx{3}{4}$ & $\CCx{4}{4}$
                \\\hline
                $\dim$ & 0 & 0& 1 & 3
        \end{tabular} 
\end{center}
The Euler characteristic number of $\ds \HGFs{\bullet}{2} = 1$ and  
the Euler characteristic number of $\ds \HGFs{\bullet}{4} = 2$.   
Those results are the first part of Theorem 1. 

\medskip

When weight=6, we again 
apply Proposition \ref{prop:lambda:s} to our decompositions, we see that 
\begin{align*}
\dim \CCx{2}{6} &= \dim \left(\Lambda^2\frakS{5}\right)^{triv}
 \mathop{=}^{ (\ref{ni:go})} 1 \\ 
 \dim \CCx{3}{6} &= 
 \dim \left((\Lambda^2\frakS{3})\otimes \frakS{6} \right)^{triv}
 +\dim \left(\frakS{3}\otimes \frakS{4}\otimes \frakS{5} \right)^{triv}
 +\dim \left(\Lambda^3\frakS{4} \right)^{triv}
 \\&
 = \dim \left((\Lambda^2\frakS{3})^{[6,0,0]}  \right)
 +\dim \left(\frakS{3}\otimes \frakS{4}\otimes \frakS{5} \right)^{triv}
 + 0 
 = \dim \left(\frakS{3}\otimes \frakS{4}\otimes \frakS{5} \right)^{triv}
 \\
 \dim \CCx{4}{6} &= 
 \dim \left((\Lambda^3\frakS{3}) \otimes \frakS{5} \right)^{triv}
 + \dim \left((\Lambda^2\frakS{3}) \otimes 
 (\Lambda^2\frakS{4}) \right)^{triv}
 \\
 &= 
 \dim \left((\Lambda^3\frakS{3})^{[5,0,0]}  \right)
 \\& \quad 
 + \dim \left((\Lambda^2\frakS{3}) \otimes 
         (V_{[2,0,0]} + V_{[3,1,0]} + V_{[4,2,0]} + V_{[5,3,0]} + V_{[6,0,0]}
         + V_{[7,1,0]} 
 ) \right)^{triv}
 \\& 
 = 
 0 
 + \dim \left((\Lambda^2\frakS{3})^{  
 [2,0,0]}\right) 
 + \dim \left((\Lambda^2\frakS{3})^{  
[3,1,0]} \right) 
 + \dim \left((\Lambda^2\frakS{3})^{  
         [4,2,0]} \right) 
 \\&\quad 
 + \dim \left((\Lambda^2\frakS{3})^{  
         [5,3,0]} \right) 
 + \dim \left((\Lambda^2\frakS{3})^{  
         [6,0,0]}\right) 
 + \dim \left((\Lambda^2\frakS{3})^{  
         [7,1,0]}  \right)
 \\& 
 = 0 
 \end{align*} 

When $n=3$ we see the crystal base of $[3,0,0]$
 consists of the next 56 tableaux. 
 \newcommand{\myoban}[3]{\begin{tabular}{|*{3}{c|}}\hline
$ #1$ &$ #2$ &$ #3$\\\hline
 \end{tabular}}
 \begin{small}
\begin{align*} 
& \myoban{1}{1}{1}, \myoban{1}{1}{2}, \myoban{1}{2}{2}, \myoban{2}{2}{2}, 
\myoban{1}{1}{3}, \myoban{1}{2}{3}, \myoban{2}{2}{3}, \myoban{1}{3}{3}, 
\\& 
\myoban{2}{3}{3},\myoban{3}{3}{3},\myoban{1}{1}{\bar{3}},\myoban{1}{2}{\bar{3}}, 
\myoban{2}{2}{\bar{3}},\myoban{1}{3}{\bar{3}},\myoban{2}{3}{\bar{3}},\myoban{3}{3}{\bar{3}},
\\&
\myoban{1}{\bar{3}}{\bar{3}}, \myoban{2}{\bar{3}}{\bar{3}}, 
\myoban{3}{\bar{3}}{\bar{3}}, \myoban{\bar{3}}{\bar{3}}{\bar{3}}, 
\myoban{1}{1}{\bar{2}}, \myoban{1}{2}{\bar{2}}, \myoban{2}{2}{\bar{2}}, 
\myoban{1}{3}{\bar{2}},\\&
\myoban{2}{3}{\bar{2}}, \myoban{3}{3}{\bar{2}}, 
\myoban{1}{\bar{3}i}{\bar{2}}, \myoban{2}{\bar{3}}{\bar{2}}, 
\myoban{3}{\bar{3}}{\bar{2}}, \myoban{\bar{3}}{\bar{3}}{\bar{2}},
\myoban{1}{\bar{2}}{\bar{2}}, \myoban{2}{\bar{2}}{\bar{2}},  \\& 
\myoban{3}{\bar{2}}{\bar{2}}, \myoban{\bar{3}}{\bar{2}}{\bar{2}}, 
\myoban{\bar{2}}{\bar{2}}{\bar{2}}, \myoban{1}{1}{\bar{1}}, 
\myoban{1}{2}{\bar{1}}, \myoban{2}{2}{\bar{1}}, \myoban{1}{3}{\bar{1}}, 
\myoban{2}{3}{\bar{1}}, \\&
\myoban{3}{ 3}{ \bar{1}}, \myoban{1}{ \bar{3}}{ \bar{1}},
\myoban{ 2}{ \bar{3}}{ \bar{1}}, \myoban{3}{ \bar{3}}{ \bar{1}},
\myoban{\bar{3}}{ \bar{3}}{ \bar{1}}, \myoban{1}{\bar{2}}{\bar{1}}, 
\myoban{2}{\bar{2}}{\bar{1}}, \myoban{3}{\bar{2}}{\bar{1}}, \\&
\myoban{\bar{3}}{\bar{2}}{\bar{1}}, \myoban{\bar{2}}{\bar{2}}{\bar{1}}, 
\myoban{1}{\bar{1}}{\bar{1}}, \myoban{2}{\bar{1}}{\bar{1}}, 
\myoban{3}{\bar{1}}{\bar{1}}, \myoban{\bar{3}}{\bar{1}}{\bar{1}}, 
\myoban{\bar{2}}{\bar{1}}{\bar{1}}, \myoban{\bar{1}}{\bar{1}}{\bar{1}} 
\end{align*}\end{small}

\newcommand{\tume}{\phantom{$\bar{2}$}}

By the Crystal Base theory, we get the irreducible decomposition of 
$\ds \frakS{3}\otimes \frakS{4} \cong
\frakS{4}\otimes \frakS{3}$  by the action of each tableau to the Young
diagram 
 \begin{tabular}{|*{4}{c|}}\hline & & & \\\hline \end{tabular}.  
The action of crystal basis of type 
 \begin{tabular}{|*{3}{c|}}\hline \tume & \tume & 3 \\\hline \end{tabular},     
 \begin{tabular}{|*{3}{c|}}\hline\tume &\tume & $\bar{3}$ \\\hline \end{tabular},    
 \begin{tabular}{|*{3}{c|}}\hline\tume &\tume & $\bar{2}$ \\\hline \end{tabular},    
 \begin{tabular}{|*{3}{c|}}\hline\tume & 3 & $\bar{1}$ \\\hline \end{tabular},    
 \begin{tabular}{|*{3}{c|}}\hline\tume & $\bar{2}$ & $\bar{1}$ \\\hline
 \end{tabular} or     
 \begin{tabular}{|*{3}{c|}}\hline\tume & $\bar{3}$ & $\bar{1}$ \\\hline
 \end{tabular} for     
 \begin{tabular}{|*{4}{c|}}\hline & & & \\\hline \end{tabular} does 
 provide nothing.  Avoiding those actions, we have 
 \begin{kmProp}
\begin{align*}
        \frakS{3}\otimes \frakS{4} \cong &
 [4,0,0] \cdot \myoban{1}{1}{1} + 
  [4,0,0] \cdot \myoban{1}{1}{2} + 
  [4,0,0] \cdot \myoban{1}{2}{2} + 
  [4,0,0] \cdot \myoban{2}{2}{2} \\&
  + [4,0,0] \cdot \myoban{1}{1}{\bar{1}} + 
  [4,0,0] \cdot \myoban{1}{2}{\bar{1}} + 
  [4,0,0] \cdot \myoban{2}{2}{\bar{1}} + 
  [4,0,0] \cdot \myoban{1}{\bar{1}}{\bar{1}} \\&
  + [4,0,0] \cdot \myoban{2}{\bar{1}}{\bar{1}} + 
  [4,0,0] \cdot \myoban{\bar{1}}{\bar{1}}{\bar{1}}  
\\ \cong &
V_{[7,0,0]} +
V_{[6,1,0]} +
V_{[5,2,0]} +
V_{[4,3,0]} +
V_{[5,0,0]} +
V_{[4,1,0]} +
V_{[3,2,0]} +
V_{[3,0,0]} +
V_{[2,1,0]} + 
V_{[1,0,0]} 
\\ \cong &
V_{[1,0,0]} +
V_{[2,1,0]} +
V_{[3,0,0]} +
V_{[3,2,0]} +
V_{[4,1,0]} +
V_{[4,3,0]} +
V_{[5,0,0]} +
V_{[5,2,0]} +
V_{[6,1,0]} +
V_{[7,0,0]} 
\end{align*}
\end{kmProp}
And so
$$ 
 \dim \CCx{3}{6} = 
 \dim \left(\frakS{3}\otimes \frakS{4}\otimes \frakS{5} \right)^{triv}
 = \dim \left( \frakS{3}\otimes \frakS{4} \right)^{[5,0,0]} = 1 \ .
$$

So far, 
concerning to $\ds \dim \CCx{\bullet}{6} $, we have got the dimensions for degree
less than 6 as follows:   
\begin{center}
        \begin{tabular}{c|*{5}c}
                & $\ds\CCx{1}{6}$ & $\CCx{2}{6}$ 
                & $\ds\CCx{3}{6}$ & $\CCx{4}{6}$
                & $\ds\CCx{5}{6}$ 
                \\\hline
                $\dim$ & 0 & 1& 1 & 0& 4    
        \end{tabular}
\end{center}

By Proposition \ref{prop:lambda:s}, 
$\ds\CCx{6}{6} = \left(\Lambda^6\frakS{3}\right)^{triv }$ and  
$\ds\dim \left(\Lambda^6\frakS{3}\right) = 32468436$.  The number of possibilities of
the maximal weight vectors is 146.  The subspace of 
$\ds\Lambda^6\frakS{3}$ generated by $\ds 
z_{A_1} \wedge z_{A_2} \wedge z_{A_3} \wedge z_{A_4} \wedge z_{A_5} \wedge 
z_{A_6}$ with 
$\ds 
{A_1} < {A_2} < {A_3} <  {A_4} < {A_5} < {A_6}$, where $\ds A_j$ are 6-dim 
vector of non-negative integer with the sum of them is 3, and
$\ds 
\sum_{j=1}^6 (A_j(1) - A_j(6)) =  
\sum_{j=1}^6 (A_j(2) - A_j(5)) =  
\sum_{j=1}^6 (A_j(3) - A_j(4)) =  0$ 
has the dimension
$204894$. By the condition that a general element is invariant under 
$U$, we have $1553660$ linear equations
with  $204894$ unknown variables.   Solve this huge linear equations with 
help of symbol calculus, we see the dimension of the null space (the kernel
space) is $7$.  
Thus, we have 
the table of dimensions of each cochain complexes with weight 6 as
\begin{center}
        \begin{tabular}{c|*{6}c}
                & $\ds\CCx{1}{6}$ & $\CCx{2}{6}$ 
                & $\ds\CCx{3}{6}$ & $\CCx{4}{6}$
                & $\ds\CCx{5}{6}$ & $\CCx{6}{6}$
                \\\hline
                $\dim$ & 0 & 1& 1 & 0& 4 & \textbf{7}  
        \end{tabular}
\end{center}
and 
the Euler characteristic number of $\ds \HGFs{\bullet}{6} = 3$.   
Those are the second part of Theorem 1.  
In fact, we have the complete irreducible decomposition of 
$\ds\Lambda^{6}\frakS{3}$ as next.  
\begin{kmProp}\label{roku:san}
\begin{align*}
\Lambda^6\frakS{3} \cong& 
V_{[12,2,0]} + 
V_{[12,0,0]} + 
V_{[11,4,1]} + 
V_{[11,3,2]} + 
3 V_{[11,3,0]} + 
2 V_{[11,2,1]} + 
V_{[11,1,0]} + 
V_{[10,6,0]} + 
2 V_{[10,5,1]} \\&  + 
V_{[10,4,4]} + 
3 V_{[10,4,2]} + 
4 V_{[10,4,0]} + 
2 V_{[10,3,3]} + 
7 V_{[10,3,1]} + 
4 V_{[10,2,2]} + 
7 V_{[10,2,0]} + 
5 V_{[10,1,1]} \\& + 
3 V_{[10,0,0]}  + 
V_{[9,7,2]} + 
V_{[9,6,3]} + 
3 V_{[9,6,1]} + 
4 V_{[9,5,2]} + 
8 V_{[9,5,0]} + 
5 V_{[9,4,3]} + 
14 V_{[9,4,1]}\\& + 
14 V_{[9,3,2]} + 
13 V_{[9,3,0]} + 
18 V_{[9,2,1]} + 
14 V_{[9,1,0]} + 
V_{[8,8,0]} + 
3 V_{[8,7,1]} + 
V_{[8,6,4]} + 
6 V_{[8,6,2]}\\& + 
4 V_{[8,6,0]} + 
6 V_{[8,5,3]} + 
16 V_{[8,5,1]} + 
3 V_{[8,4,4]} + 
22 V_{[8,4,2]} + 
24 V_{[8,4,0]} + 
10 V_{[8,3,3]} + 
39 V_{[8,3,1]}\\& + 
20 V_{[8,2,2]} + 
30 V_{[8,2,0]} + 
20 V_{[8,1,1]} + 
9 V_{[8,0,0]} + 
V_{[7,7,4]} + 
V_{[7,7,2]} + 
5 V_{[7,7,0]} + 
5 V_{[7,6,3]} + 
12 V_{[7,6,1]}\\& + 
3 V_{[7,5,4]} + 
24 V_{[7,5,2]} + 
15 V_{[7,5,0]} + 
19 V_{[7,4,3]} + 
49 V_{[7,4,1]} + 
41 V_{[7,3,2]} + 
55 V_{[7,3,0]} + 
57 V_{[7,2,1]}\\& + 
29 V_{[7,1,0]} + 
3 V_{[6,6,4]} + 
7 V_{[6,6,2]} + 
15 V_{[6,6,0]} + 
V_{[6,5,5]}  + 
14 V_{[6,5,3]} + 
33 V_{[6,5,1]} + 
7 V_{[6,4,4]} \\& + 
53 V_{[6,4,2]} + 
41 V_{[6,4,0]} + 
26 V_{[6,3,3]} + 
82 V_{[6,3,1]} + 
41 V_{[6,2,2]} + 
68 V_{[6,2,0]} + 
41 V_{[6,1,1]} + 
15 V_{[6,0,0]}\\& + 
5 V_{[5,5,4]} + 
16 V_{[5,5,2]} + 
28 V_{[5,5,0]} + 
24 V_{[5,4,3]} + 
      62  V_{[5, 4, 1]} + 
      70  V_{[5, 3, 2]} + 
      56  V_{[5, 3, 0]} + 
      84  V_{[5, 2, 1]}\\& + 
      57  V_{[5, 1, 0]} + 
       7  V_{[4, 4, 4]} + 
      25  V_{[4, 4, 2]} + 
      40  V_{[4, 4, 0]} + 
      22  V_{[4, 3, 3]} + 
      73  V_{[4, 3, 1]} + 
      48  V_{[4, 2, 2]} + 
      54  V_{[4, 2, 0]} \\& + 
      46  V_{[4, 1, 1]} + 
      25  V_{[4, 0, 0]} + 
      21  V_{[3, 3, 2]} + 
      46  V_{[3, 3, 0]} + 
      51  V_{[3, 2, 1]} + 
      27  V_{[3, 1, 0]} + 
       6  V_{[2, 2, 2]} + 
       32  V_{[2, 2, 0]}\\& + 
      15  V_{[2, 1, 1]} + 
       3  V_{[2, 0, 0]} + 
       14  V_{[1, 1, 0]} + 
       7  V_{[0, 0, 0]} 
\end{align*} 
\end{kmProp}

\section{Betti numbers}
The definition of Betti number is 
$\ds  b^{\bullet}_{w} := \dim 
\HGFs{\bullet}{w}$ and so \begin{align*}
 b^{j}_{w} 
 =& \dim\left( \ker ( d_j : \CCx{j}{w}\rightarrow \CCx{j+1}{w} )\right) 
- \dim\left( d_{j-1} ( \CCx{j-1}{w})  \right)
\\
=& 
\dim\left( \CCx{j}{w} \right) 
- \dim\left( d_{j} ( \CCx{j}{w})  \right)
- \dim\left( d_{j-1} ( \CCx{j-1}{w})  \right)
\\
=& 
\dim\left( \CCx{j}{w} \right) 
- \rank(  d_{j}) - \rank(  d_{j-1})
\end{align*}
From the definition of the coboundary operator, the operator has
skew-derivation of degree 1 and 
for each 1-cochain $\sigma$,  $\ds \myd( \sigma)(f,g) := - \langle
\sigma, \Pkt{f}{g}\rangle$ where $\ds \Pkt{f}{g}$ is the Poisson bracket of
the standard symplectic structure of $\ds\mR^{6}$.  

When the weight is 2, we see $\ds b_{2}^{j} = \dim 
\CCx{j}{2}$ directly.    

In general, to investigate the Betti numbers, we shall know the rank of the
consecutive coboundary operators. For that purpose, we have to prepare concrete
bases of the cochain complexes.

As commented briefly just before Proposition \ref{roku:san},  
each cochain is generated by $\ds z_A $ where $A$ is 6-dimensional 
multi-index with non-negative integers. We may regard $\ds z_A$ as the dual 
of $\ds 
\frac{x_1{}^{a_1}} { a_1 !} 
\frac{x_2{}^{a_2}} { a_2 !} 
\cdots 
\frac{x_6{}^{a_6}} { a_6 !} 
$ where $\ds A = (a_1, a_2, a_3, a_4, a_5, a_6)$ with 
$\ds a_j \in \mN^{+} $ ($j=1..6$).     
$\ds z_A \in\frakS{h}$ is equivalent to 
$\ds \sum_{j=1}^6 a_j = h$.  

In this paper, we use the abbreviation $\ds Z^{a_1 a_2 a_3}_{a_4 a_5 a_6}$
instead of $\ds z[ a_1, a_2, a_3, a_4, a_5, a_6]$.    

\newcommand{\km}{\wedge}
\newcommand{\W}[2]{Z^{#1}_{#2}}

When the weight is 4, 
as a basis of $\ds \CCx{3}{4}$ 
we have 
\begin{small}
\begin{align*} & 
-\W{100}{101}\km\W{000}{102}\km\W{202}{000}
-\frac{1}{2} \W{001}{011}\km\W{000}{102}\km\W{310}{000}
+\frac{1}{2} \W{200}{001}\km\W{000}{003}\km\W{300}{001}
+\frac{1}{2} \W{110}{001}\km\W{000}{003}\km\W{300}{010} \\&
+\frac{1}{2} \W{110}{010}\km\W{000}{003}\km\W{210}{010}
+\frac{1}{2} \W{011}{001}\km\W{000}{102}\km\W{300}{010}
+\frac{1}{2} \W{002}{001}\km\W{000}{102}\km\W{300}{100}
+\frac{1}{2} \W{200}{010}\km\W{000}{003}\km\W{210}{001} \\&
+ (1228 \text{ terms}) \\&
-\frac{1}{2} \W{002}{100}\km\W{101}{010}\km\W{010}{201}
+\frac{1}{2} \W{200}{100}\km\W{020}{001}\km\W{001}{021}
-\W{102}{000}\km\W{000}{111}\km\W{110}{101}
-\frac{1}{2} \W{000}{120}\km\W{101}{001}\km\W{120}{001} \\&
-\W{110}{100}\km\W{101}{100}\km\W{001}{012}
+\W{010}{110}\km\W{011}{001}\km\W{101}{110} 
\end{align*} 
\end{small} 
(You may see the full form of the above vector at 

\texttt{http://www.math.akita-u.ac.jp/\~{}mikami/GKF\_R6/c3w4\_d6.pdf}, where 
we denote  $\ds\W{100}{101}\km\W{000}{102}\km\W{202}{000}$ by 
\text{w(100101,000102,202000)}) for example,  
and the $d$-image is the next form, which is not 0. 
\begin{small}
\begin{align*}&
-6 \W{000}{120}\km\W{002}{001}\km\W{110}{001}\km\W{110}{100}
-\W{111}{000}\km\W{000}{120}\km\W{120}{000}\km\W{000}{012}
+8 \W{111}{000}\km\W{000}{120}\km\W{110}{001}\km\W{010}{011}\\&
+8 \W{111}{000}\km\W{000}{120}\km\W{110}{010}\km\W{010}{002} 
-\W{111}{000}\km\W{000}{120}\km\W{100}{011}\km\W{020}{001}
+7 \W{111}{000}\km\W{000}{120}\km\W{020}{001}\km\W{010}{020} \\&
+ (2645 \text{terms}) \\&
+\frac{1}{2} \W{120}{000}\km\W{000}{201}\km\W{010}{020}\km\W{002}{010}
+\frac{1}{2} \W{021}{000}\km\W{000}{300}\km\W{010}{020}\km\W{002}{010}
-2 \W{010}{020}\km\W{110}{001}\km\W{200}{001}\km\W{010}{011} 
\end{align*} 
\end{small}        
(cf.\ \texttt{http://www.math.akita-u.ac.jp/\~{}mikami/GKF\_R6/d-image-c3w4\_d6.pdf}
(21 pages)).    
Thus, 
$\ds\rank( d_3 )=1$ and this implies 
$\ds b^{3}_{4} =0$ and 
$\ds b^{4}_{4} =2$ and the others are 0.   This completes a proof to 
the first half of Theorem 2.  

When the weight is 6, 
as a basis of $\ds \CCx{2}{6}$ 
we have 
\begin{align*} &
-\frac{1}{2} \W{121}{010}\km\W{010}{121}
-\frac{1}{2} \W{121}{001}\km\W{100}{121}
+\frac{1}{4} \W{120}{200}\km\W{002}{021}
+\frac{1}{2} \W{120}{110}\km\W{011}{021} 
+\frac{1}{2} \W{120}{101}\km\W{101}{021} \\&
+\frac{1}{4} \W{120}{020}\km\W{020}{021}
+\frac{1}{2} \W{120}{011}\km\W{110}{021}
+\frac{1}{6} \W{113}{000}\km\W{000}{311} 
-\frac{1}{2} \W{112}{100}\km\W{001}{211}
-\frac{1}{2} \W{112}{010}\km\W{010}{211}\\&
+ (111 \text{ terms}) \\&
-\frac{1}{6} \W{130}{100}\km\W{001}{031} 
-\frac{1}{6} \W{130}{010}\km\W{010}{031}
-\frac{1}{6} \W{130}{001}\km\W{100}{031}
+\frac{1}{4} \W{122}{000}\km\W{000}{221}
-\frac{1}{2} \W{121}{100}\km\W{001}{121}
\end{align*} 
(cf.\ \texttt{http://www.math.akita-u.ac.jp/\~{}mikami/GKF\_R6/c2w6\_d6.pdf})  
and the $d$-image is the next form, which is not 0. 
\begin{small}
\begin{align*}&
-\frac{1}{4} \W{200}{030}\km\W{030}{000}\km\W{010}{012}
-\frac{1}{2} \W{200}{030}\km\W{030}{001}\km\W{010}{011}
+\frac{1}{12} \W{200}{030}\km\W{030}{010}\km\W{010}{002}
-\frac{1}{12} \W{200}{030}\km\W{040}{000}\km\W{000}{012} \\&
+\frac{1}{12} \W{200}{030}\km\W{000}{102}\km\W{031}{000}
-\frac{1}{12} \W{200}{030}\km\W{001}{002}\km\W{030}{100}
+\frac{1}{2} \W{200}{030}\km\W{010}{101}\km\W{021}{001}
+\frac{1}{4} \W{200}{030}\km\W{010}{102}\km\W{021}{000} \\& 
+( 6510 \text{ terms}) \\&
-\frac{1}{4} \W{201}{020}\km\W{020}{110}\km\W{010}{002}
-\frac{1}{4} \W{201}{020}\km\W{020}{010}\km\W{010}{102}
-\frac{1}{2} \W{201}{020}\km\W{020}{011}\km\W{010}{101}
+\frac{1}{2} \W{201}{020}\km\W{020}{100}\km\W{010}{012}
\end{align*} 
\end{small}
(cf.\ 
\texttt{http://www.math.akita-u.ac.jp/\~{}mikami/GKF\_R6/d-image-c2w6\_d6.pdf}
(41 pages)).   

Thus, we see that $\ds b^{1}_{6} = b^{2}_{6} = b^{3}_{6} = b^{4}_{6} = 0$.  

In order to know the rank of $\ds d : \CCx{5}{6}\longrightarrow \CCx{6}{6}$,
we have to fix some bases of both spaces, whose dimensions are 4 and 7
respectively.    
Concerning to finding a basis of $\ds \CCx{6}{6} = \left(
\Lambda^{6}\frakS{3}\right)^{triv}$, we can use our strategy of splitting 
$\ds \Lambda^{p}\frakS{q}$ and we have got them, say $\ds    
  \mathbf{r}_1, 
   \mathbf{r}_2, 
   \mathbf{r}_3, 
   \mathbf{r}_4, 
   \mathbf{r}_5, 
   \mathbf{r}_6, 
   \mathbf{r}_7$.  

\begin{kmProp}\label{P:6}
The number of summands of 
 $\ds \mathbf{r}_1$ is 150340. Here we may write  
 $\ds \#\mathbf{r}_1=150340$. Then we have  
 $\ds \#\mathbf{r}_2=21612$, 
 $\ds \#\mathbf{r}_3=153466$, 
 $\ds \#\mathbf{r}_4=148660$, 
 $\ds \#\mathbf{r}_5=155512$, 
 $\ds \#\mathbf{r}_6=3276$, and 
 $\ds \#\mathbf{r}_7=148600$.  
\end{kmProp}
You will understand it is hard to show them in this paper. 
You will see them at\\ 
\texttt{http://www.math.akita-u.ac.jp/\~{}mikami/GKF\_R6/c6w6-j\_d6.pdf}
(j=1,\ldots, 7).  
The first vector spends 1567 pages to show, 
the second vector\kmcomment{spends}  226 pages, 
the third vector\kmcomment{spends} 1599 pages, 
the fourth vector\kmcomment{spends} 1549 pages, 
the fifth vector\kmcomment{spends} 1620 pages, 
the sixth vector\kmcomment{spends}   35 pages, and 
the last seventh vector spends   1548 pages.  

\bigskip

If we try to 
find a concrete basis of 
$\ds\CCx{5}{6} =\left((\Lambda^{4}\frakS{3})\otimes\frakS{4}\right)^{triv}$ 
directly, then 
from $\ds \dim 
\left((\Lambda^{4}\frakS{3})\otimes\frakS{4}\right)
= 
{\dim \frakS{3} \choose 4 }\cdot \dim \frakS{4}$, 
we have to handle 46278540  variables. Bus,   
as we discussed in the section of Crystal base, trivial representation spaces
live in $\ds \left(\Lambda^{4}\frakS{3}\right)^{[4,0,0]} \otimes\frakS{4}$.
Then we only deal with $\ds 126^2 = 15876$ variables ($\ds\dim
\frakS{4} = 126$).  
Proposition \ref{prop:lambda:s} (\ref{yon:san}) tells us that the
multiplicity is 4. Thus, we first fix the 4 linearly independent maximal 
vectors, say $\ds w_1, w_2, w_3, w_4$  of weight $[4,0,0]$ of 
$\ds \Lambda^{4}\frakS{3}$.  
Let us say $\ds W_{j}$ be the irreducible subspace of the maximal vector $\ds
w_j$ ($j=1,2,3,4$).  Using the Fact (2), 
We get a basis of the space $\ds W_{j}$ of 126-dimensional.  
On each the tensor product space 
$\ds  W_j \otimes \frakS{4} $,  
after a long calculation, we get the maximal vector of weight $[0,0,0]$. 
Those are 
member of 
a concrete basis of 
$\ds \CCx{5}{6} = \left((\Lambda ^{4}\frakS{3} )\wedge
\frakS{4}\right)^{triv}$. The next are small part of them. 

\begin{kmProp}\label{P:5}
\begin{small} 
\begin{align*} &
-6 \W{210}{000}\km\W{002}{010}\km\W{000}{120}\km\W{000}{012}\km\W{030}{100}
+12 \W{210}{000}\km\W{200}{010}\km\W{100}{011}\km\W{000}{012}\km\W{020}{002} 
\\&
+ ( 32141\text{ terms}) \\& 
+12 \W{002}{001}\km\W{100}{101}\km\W{110}{010}\km\W{010}{011}\km\W{101}{200}
+12 \W{210}{000}\km\W{200}{010}\km\W{010}{011}\km\W{000}{012}\km\W{020}{011} 
\end{align*}
\begin{align*}&
-12 \W{110}{001}\km\W{011}{100}\km\W{010}{110}\km\W{100}{002}\km\W{101}{020}
-6 \W{002}{100}\km\W{020}{100}\km\W{010}{011}\km\W{010}{110}\km\W{101}{020}
\\&
+( 177494 \text{ terms}) \\& 
+6 \W{120}{000}\km\W{012}{000}\km\W{100}{011}\km\W{000}{120}\km\W{100}{102}
+12 \W{100}{101}\km\W{002}{100}\km\W{110}{010}\km\W{011}{010}\km\W{100}{102} 
\end{align*}
\begin{align*}&
20 \W{200}{010}\km\W{102}{000}\km\W{000}{201}\km\W{100}{011}\km\W{020}{002}
+6 \W{111}{000}\km\W{100}{101}\km\W{100}{110}\km\W{000}{003}\km\W{201}{001} 
\\&
+ ( 203321\text{ terms}) \\& 
-3 \W{011}{010}\km\W{020}{100}\km\W{200}{001}\km\W{000}{012}\km\W{110}{020}
+2 \W{101}{100}\km\W{001}{200}\km\W{010}{200}\km\W{100}{011}\km\W{003}{001} 
\end{align*}
\begin{align*}&
-6 \W{030}{000}\km\W{101}{100}\km\W{200}{010}\km\W{000}{111}\km\W{001}{012}
-3 \W{201}{000}\km\W{000}{300}\km\W{020}{010}\km\W{010}{002}\km\W{002}{020} 
\\&
+ ( 188471\text{ terms}) \\& 
-6 \W{011}{010}\km\W{000}{102}\km\W{010}{110}\km\W{001}{020}\km\W{220}{000}
-6 \W{300}{000}\km\W{100}{101}\km\W{001}{110}\km\W{001}{002}\km\W{110}{002}
\end{align*} 
\end{small} 
(cf. \texttt{http://www.math.akita-u.ac.jp/\~{}mikami/GKF\_R6/c5w6-j\_d6.pdf}
(j=1,2,3,4))  
\end{kmProp}

\begin{thm}\rm 
        Let $\ds \mathbf{w}_j$ ($j=1,\ldots,4$) be the basis described in 
        Proposition \ref{P:5} and  
        $\ds \mathbf{r}_k$ ($k=1,\ldots,7$) be the basis in 
        Proposition \ref{P:6}.   
        Then the coboundary operator $\myd$ has the representation below:
$$
[ \myd( \mathbf{w}_1), 
        \myd( \mathbf{w}_2), 
        \myd( \mathbf{w}_3), 
\myd( \mathbf{w}_4)] = 
[  \mathbf{r}_1, 
   \mathbf{r}_2, 
   \mathbf{r}_3, 
   \mathbf{r}_4, 
   \mathbf{r}_5, 
   \mathbf{r}_6, 
   \mathbf{r}_7] 
  \begin{bmatrix} 
72 &     0 &        -864 &    -576 \\       
-384 &   1152 &     -11520 &  -5376 \\
0 &      -48 &      384 &     -192 \\
192 &    -192 &     1920 &    -576 \\
0 &      -168 &     2400 &    1824 \\ 
-276 &   720 &      -11712 &  -2976 \\
-84  &   180 &      -1728 &   -816 \\ 
\end{bmatrix}
$$
And so we have $\ds\rank(\myd) = 2$. 
\end{thm}\rm 



\def\cprime{$'$} \def\cprime{$'$}

}

\begin{thebibliography}{10}

\bibitem{fulton:harris}
William Fulton and Joe Harris.
\newblock {\em Representation {T}heory}, volume 129 of {\em Graduate Texts in
  Math}.
\newblock Springer Verlag, 1991.

\bibitem{MR0312531}
I.~M. Gel{\cprime}fand, D.~I. Kalinin, and D.~B. Fuks.
\newblock The cohomology of the {L}ie algebra of {H}amiltonian formal vector
  fields.
\newblock {\em Funkcional. Anal. i Prilo\v zen.}, 6(3):25--29, 1972.

\bibitem{goodman:wallach}
Roe Goodman and R.~Nolan Wallach.
\newblock {\em Symmetry, {R}epresentations, and {I}nvariants}, volume 255 of
  {\em Graduate Texts in Math}.
\newblock Springer Verlag, 2009.

\bibitem{kashi:french}
Masaki Kashiwara.
\newblock {\em Bases {C}ristallines des {G}roupes {Q}uantiques}.
\newblock Soci\'et\'e {M}ath\'ematique de {F}rance, Paris, 2002.
\newblock Cours {S}p\'ecialis\'es 9.

\bibitem{kashi:nakashima}
Masaki Kashiwara and Toshiki Nakashima.
\newblock Crystal {G}raphs for {R}epresentations of the q-{A}nalogue of
  {C}lassical {L}ie {A}lgebras.
\newblock {\em J. of Algebra}, 165:295--345, 1994.

\bibitem{Kont:RW:MR1671725}
Maxim Kontsevich.
\newblock Rozansky-{W}itten invariants via formal geometry.
\newblock {\em Compositio Math.}, 115(1):115--127, 1999.

\bibitem{KOT:MORITA}
D.~Kotschick and S.~Morita.
\newblock The {G}el'fand-{K}alinin-{F}uks class and characteristic classes of
  transversely symplectic foliations.
\newblock {\em arXiv:0910.3414}, October 2009.

\bibitem{MR1354144}
I.~G. Macdonald.
\newblock {\em Symmetric functions and {H}all polynomials}.
\newblock Oxford Mathematical Monographs. The Clarendon Press Oxford University
  Press, New York, second edition, 1995.
\newblock With contributions by A. Zelevinsky, Oxford Science Publications.

\bibitem{metoki:shinya}
Shinya Metoki.
\newblock {\em Non-trivial cohomology classes of Lie algebras of volume
  preserving formal vector fields}.
\newblock PhD thesis, Univ.\ of Tokyo, 2000.

\bibitem{Mik:Nak}
K.~Mikami and Y.~Nakae.
\newblock {L}ower weight {G}el'fand-{K}alinin-{F}uks cohomology groups of the
  formal {H}amiltonian vector fields on ${\mR}^4$.
\newblock {\em J.\ Math. Sci. Univ. Tokyo}, 19:1--18, 2012.

\bibitem{M:N:K}
K.~Mikami, Y.~Nakae, and H.~Kodama.
\newblock {H}igher weight {G}el'fand-{K}alinin-{F}uks classes of formal
  {H}amiltonian vector fields of symplectic ${\mR}^2$.
\newblock {\em arXiv:1210.1662v1}, October 2012.

\bibitem{KM:affirm}
Kentaro Mikami.
\newblock {A}n affirmative answer to a conjecture for {M}etoki class.
\newblock preprint, October 2012.

\bibitem{morita:text:eng}
Shigeyuki Morita.
\newblock {\em Geometry of characteristic classes}, volume 199 of {\em
  Translations of Math. Monograph}.
\newblock AMS, 2001.

\bibitem{morita:text}
Shigeyuki Morita.
\newblock {\em Characteristic classes and Geometry}.
\newblock Iwanami, 2008.
\newblock (in Japanase).

\bibitem{nakashima}
Toshiki Nakashima.
\newblock Crystal {B}ase and a {G}eneralization of the
  {L}ittlewood-{R}ichardson {R}ule for the {C}lassical {L}ie {A}lgebras.
\newblock {\em Commun. Math. Phys.}, 154:215--243, 1993.

\bibitem{okada:text}
Soichi Okada.
\newblock {\em Representation theory of classical groups and Combinatorics}.
\newblock Baifukan, 2006.
\newblock (in Japanase).

\bibitem{M:Takamura}
Masashi Takamura.
\newblock The relative cohomology of formal contact vector fields with respect
  to formal {P}oisson vector fields.
\newblock {\em J.~Math.~Soc.~Japan}, 60(1):117--125, 2008.

\end{thebibliography}
\end{document}